\documentclass[a4paper, 10pt]{article}
\usepackage{graphicx}

\usepackage{alltt}

\usepackage{amsmath,amssymb}

\usepackage{yhmath} 
\usepackage{srctex} 
\usepackage[dvipsnames]{color}

\usepackage{srctex} 

\def\m#1{\mathsf{#1}} 

\newcommand{\cl}[2]{\ensuremath{\mathit{Cl}_{#1,#2}}}

\DeclareMathOperator{\Det}{Det} 
\newcommand{\Adj}{{\mathop{\textrm{Adj}}}}

\newcommand{\bbR}{\ensuremath{\mathbb{R}}}
\newcommand{\bbC}{\ensuremath{\mathbb{C}}}

\newcommand{\reverse}[1]{\widetilde{#1}}
\newcommand{\gradeinverse}[1]{\wideparen{#1}}

\newcommand{\ii}{\mathrm{i}}

\newcommand{\ee}{\mathrm{e}} 
\newcommand{\dd}{\mathrm{d}} 

\def\A{\mathsf{A}}

\def\m#1{\mathsf{#1}}

\def\e#1{\mathbf{e}_{#1}} 

\newcommand{\ba}{\ensuremath{\mathbf{a}}}
\newcommand{\bb}{\ensuremath{\mathbf{b}}}
\newcommand{\bB}{\ensuremath{\mathbf{B}}}

\newcommand{\cA}{\ensuremath{\mathcal{A}}}
\newcommand{\cB}{\ensuremath{\mathcal{B}}}

\renewcommand{\d}{\ensuremath{\mathbin{\cdot}}} 
\newcommand{\w}{\ensuremath{\mathbin{\wedge}}} 

\usepackage{color}

\newtheorem{thm}{Theorem}[section]

 \numberwithin{equation}{section}

\begin{document}

\begin{center}
Nonlinear Analysis: Modelling and Control, Vol. vv, No. nn, YYYY\\
\copyright\ Vilnius University\\[24pt]
\LARGE \textbf{Exponentials of general multivector (MV)\\ in 3D
Clifford algebras}\\[6pt]
\small
\textbf{Adolfas Dargys$^{*}$, Art{\=u}ras Acus$^{**}$}\\[6pt]
$^{*}$Center for Physical Sciences and Technology, Semiconductor
Physics Institute, \\ Saul{\.e}tekio 3, LT-10257
Vilnius, Lithuania\\
adolfas.dargys@ftmc.lt\\
$^{**}$Institute of Theoretical Physics and Astronomy, Vilnius
University,\\ Saul{\.e}tekio 3, LT-10257 Vilnius, Lithuania\\
arturas.acus@tfai.vu.lt\\[6pt]

Received: date\quad/\quad Revised: date\quad/\quad Published
online: data
\end{center}



\begin{abstract}
Closed form expressions  to calculate the exponential of a general
multivector (MV)
in Clifford geometric algebras (GAs) $\cl{p}{q}$ are presented for
$n=p+q= 3$. The obtained  exponential formulas were applied to
find exact  GA trigonometric and hyperbolic functions of MV
argument. We have verified that the presented  exact formulas are
in accord with series  expansion of MV hyperbolic and
trigonometric functions. The exponentials may be applied ro solve
GA differential equations, in signal and image processing,
automatic control and robotics.
\end{abstract}

computer-aided theory.


\section{Introduction}
\label{sec:intro}

In Clifford geometric algebra (GA), the exponential functions with
the exponent represented by a simple blade are well-known and used
widely. In case of complex algebra  (the complex number algebra is
isomorphic to  \cl{0}{1} GA) the exponential can be expanded into
a trigonometric function sum  by de~Moivre's theorem. In 2D vector
space, including Hamilton quaternions, the exponential is similar
to de Moivre's formula multiplied by exponential  of the scalar
part~\cite{Gurlebeck1997,Lounesto97,Chappell2015,Josipovic2019}.
In 3D vector spaces only special cases are known. Particularly,
when the square of the blade is equal to $\pm 1$, the exponential
can be expanded in de~Moivre-type sum of trigonometric or
hyperbolic functions, respectively. However, general expansion in
a symbolic form in case of 3D algebras \cl{3}{0}, \cl{1}{2},
\cl{2}{1} and \cl{0}{3}, when the exponent is a general
multivector (MV), is more difficult. The paper~\cite{Chappell2015}
considers general properties of functions of MV variable for
Clifford algebras $n=p+q\le 3$, including the exponential
function, for this purpose using the unique properties of a
pseudoscalar $I$ in \cl{3}{0} and \cl{1}{2} algebras.
Namely, the pseudoscalar in these algebras commutes with all
MV elements and $I^2=-1$. This allows to introduce more general
functions, in particular, the polar decomposition of all
multivectors. A different approach  to resolve the problem is to
factor, if possible, the exponential into product of simpler
exponentials, for example, in the  polar
form~\cite{Hitzer2020a,Hitzer2019a}. General bivector exponentials
in \cl{4}{1} algebra were analyzed in~\cite{Cameron2004}. In
coordinate form, the difficulty is connected with the appearance
of both trigonometric and hyperbolic functions simultaneously  in
the expansion of exponentials as well as the mixing of scalar
coefficients from different grades.

In this paper a different approach which presents the exponential
in  coordinates and which is more akin to construction of de
Moivre formula was applied. Namely, to solve the problem the GA
exponential function is expanded into sum of basis elements
(grades) using for this purpose the computer algebra
(\textit{Mathematica} package). Although in this way obtained
final formulas are rather cumbersome, however, their analysis
allows to identify the obstacles in constructing the GA
coordinate-free formulas. In the paper presented formulas
can be also applied to general purpose programming languages such
as Fortran,  C$^{++}$ or Python.

In Sec.~\ref{sec:notations} the notation is introduced. The final
exponential formulas in the coordinates are presented in
Secs.~\ref{sec:Cl03}-\ref{sec:Cl21} in a form of theorems. The
particular cases that follow from general exponential formulas are
given in Sec.~\ref{sec:particularCases}. Relations of GA
exponential to GA trigonometric and hyperbolic functions is
presented in Sec.~\ref{Relations}. Possible application of the
exponential function in solving spinorial Pauli-Schr{\"o}dinger
equation are given in Sec.~\ref{sec:applications}. In
Sec.~\ref{conclusion} we discuss further development of the
problem. In the Appendix we compare finite GA series  of
trigonometric functions with the exact formulas that follow from
exponential.

\section{Notation}
\label{sec:notations} In the inverse degree lexicographic ordering
used in this paper, the  general MV in GA space  is
expanded in the orthonormal basis
$\{1,\e{1},\e{2},\e{3},\e{12},\e{13},\e{23},\e{123}\equiv I\}$,
where $\e{i}$ are basis vectors, $\e{ij}$ are the bivectors and
$I$ is the pseudoscalar.\footnote{\label{note1}An increasing order
of digits  in basis elements is used, i.e., we write $\e{13}$
instead of $\e{31}=-\e{13}$. This convention is reflected in
opposite signs of some terms in formulas.} The number of
subscripts indicates the grade. The scalar is a grade-0 element,
the vectors $\e{i}$ are the grade-1 elements, etc. In the
orthonormalized basis the geometric products of basis vectors
satisfy the anticommutation relation,
 \begin{equation}\label{anticom}
 \e{i}\e{j}+\e{j}\e{i}=\pm 2\delta_{ij}.
 \end{equation}
For \cl{3}{0} and \cl{0}{3} algebras the squares of basis vectors,
correspondingly, are $\e{i}^2=+1$ and $\e{i}^2=-1$, where
$i=1,2,3$. For mixed signature algebras such as  \cl{2}{1} and
\cl{1}{2} we have $\e{1}^2=\e{2}^2=1$, $\e{3}^2=-1$ and
$\e{1}^2=1$, $\e{2}^2=\e{3}^2=-1$, respectively. The general MV of
real Clifford algebras $\cl{p}{q}$ for $n=p+q= 3$ can be expressed
as
\begin{equation}\begin{split}\label{mvA}
\A=&\,a_0+a_1\e{1}+a_2\e{2}+a_3\e{3}+a_{12}\e{12}+a_{23}\e{23}+a_{13}\e{13}+a_{123}I\\
\equiv&\,a_0+\ba+\cA+a_{123}I,
\end{split}\end{equation} where $a_i$, $a_{ij}$ and $a_{123}$ are the real
coefficients, and $\ba=a_1\e{1}+a_2\e{2}+a_3\e{3}$ and
$\cA=a_{12}\e{12}+a_{23}\e{23}+a_{13}\e{13}$ is, respectively, the
vector and bivector. $I$~is the pseudoscalar, $I=\e{123}$.
Similarly, the exponential $\m{B}$ will be denoted as
\begin{equation}
\begin{split}\label{mvB}
\m{B}=\ee^{\m{A}}=&\,b_0+b_1\e{1}+b_2\e{2}+ba_3\e{3}+b_{12}\e{12}+b_{23}\e{23}+ba_{13}\e{13}+a_{123}I\\
\equiv&\,b_0+\bb+\cB+b_{123}I.
\end{split}
\end{equation}

We start from the $\cl{0}{3}$ geometric algebra (GA) where the
expanded exponential in the coordinate form has the simplest MV
coefficients.

\section{MV exponential in \cl{0}{3} algebra}
\label{sec:Cl03}

\begin{thm}[Exponential function of multivector
in \cl{0}{3}] \label{exp03thm} The exponential of MV
$\m{A}=a_0+a_1\e{1}+a_2\e{2}+a_3\e{3}+a_{12}\e{12}+a_{13}\e{13}+a_{23}\e{23}+a_{123}
I$ is  another MV
 $\exp(\m{A})= b_0+b_1\e{1}+b_2\e{2}+b_3\e{3}+b_{12}\e{12}+b_{13}\e{13}+b_{23}\e{23}+b_{123}I$,
  where the real coefficients are 
 \begin{align}\label{exp03F}
&\begin{aligned}
\phantom{{}_{99}}b_{0}&=\tfrac{1}{2}\ee^{a_{0}} \bigl(\ee^{a_{123}}\cos a_{+} + \ee^{-a_{123}} \cos a_{-}\bigr),\\
b_{123}&=\tfrac{1}{2}\ee^{a_{0}} \bigl(\ee^{a_{123}}\cos a_{+} -\ee^{-a_{123}}\cos a_{-}\bigr)\,
\end{aligned}
\notag\allowdisplaybreaks \\[4pt]
&\begin{aligned}
  \phantom{{}_{99}}b_{1}&=\tfrac{1}{2}\ee^{a_{0}} \Bigl(\ee^{a_{123}}(a_{1}-a_{23})\frac{\sin a_{+}}{a_{+}} + \ee^{-a_{123}}(a_{1}+a_{23})\frac{\sin a_{-}}{a_{-}}\Bigr),\\
\phantom{{}_{99}}b_{2}&=\tfrac{1}{2}\ee^{a_{0}} \Bigl(\ee^{a_{123}}(a_{2}+a_{13})\frac{\sin a_{+}}{a_{+}} +\ee^{-a_{123}}(a_{2}-a_{13})\frac{\sin a_{-}}{a_{-}}\Bigr),\\
\phantom{{}_{99}}b_{3}&=\tfrac{1}{2}\ee^{a_{0}} \Bigl(\ee^{a_{123}}(a_{3}-a_{12})\frac{\sin a_{+}}{a_{+}} +\ee^{-a_{123}}(a_{3}+a_{12})\frac{\sin a_{-}}{a_{-}}\Bigr),
\end{aligned}\notag\allowdisplaybreaks\\[4pt]
&
\begin{aligned}
\phantom{{}_{9}}b_{12}&=\tfrac{1}{2}\ee^{a_{0}} \Bigl(-\ee^{a_{123}}(a_{3}-a_{12})\frac{\sin a_{+}}{a_{+}} +\ee^{-a_{123}}(a_{3}+a_{12})\frac{\sin a_{-}}{a_{-}}\Bigr),\\
\phantom{{}_{9}}b_{13}&=\tfrac{1}{2}\ee^{a_{0}} \Bigl(\ee^{a_{123}}(a_{2}+a_{13})\frac{\sin a_{+}}{a_{+}} -\ee^{-a_{123}}(a_{2}-a_{13})\frac{\sin a_{-}}{a_{-}}\Bigr),\\
\phantom{{}_{9}}b_{23}&=\tfrac{1}{2}\ee^{a_{0}} \Bigl(-\ee^{a_{123}}(a_{1}-a_{23})\frac{\sin a_{+}}{a_{+}} +\ee^{-a_{123}}(a_{1}+a_{23})\frac{\sin a_{-}}{a_{-}}\Bigr),
\end{aligned}\allowdisplaybreaks\\
\textrm{and where}\notag\\
  a_{+}&=\sqrt{(a_{3}-a_{12})^2+(a_{2}+a_{13})^2+(a_{1}-a_{23})^2}\,,\notag \\
  a_{-}&=\sqrt{(a_{3}+a_{12})^2+(a_{2}-a_{13})^2+(a_{1}+a_{23})^2}\,.
\end{align}
When either $a_{+}=0$ or $a_{-}=0$, or the both are equal to zero
simultaneously, the formula yields special cases considered in the
Subsec.~\ref{specialcasesCL03}.
\end{thm}

\textit{Proof}.\\ The simplest way to prove the above formula
$\exp(\m{A})$ is to check explicitly its defining property:
\begin{equation}\label{definingProperty}
\left.\frac{\partial\exp(\m{A}t)}{\partial t}\right|_{t=1} = \m{A}
\exp(\m{A})=\exp(\m{A}) \m{A}, \end{equation} where $\m{A}$
is assumed to be independent of $t$. Since we have a single MV
that always commutes with itself the multiplications from left and
right by $\m{A}$ coincide. After differentiation with
  respect to scalar parameter $t$ and then setting $t=1$ we find that in
  this way obtained result indeed is $\m{A} \exp(\m{A})$. To be sure we
also checked the Eq.~\eqref{definingProperty} by series expansions
of $\exp(\m{A}t)$ up to order 6 with symbolic coefficients and up
to order 20 with random integers using for this purpose the {\it
Mathematica} package~\cite{AcusDargys2017}.

\subsection{Special cases  of Theorem~\ref{exp03thm}}\label{specialcasesCL03}
Let  $\Det(\m{A})$ be the determinant of MV
\cite{Dadbeh2011,Hitzer2016,Acus2018,Shirokov2020a}. The
determinant of the sum of vector $\ba$ and bivector $\cA$ parts of
$\m{A}$ simplifies to
\begin{equation}
 \begin{split}\label{aPaM03}
   \Det(\ba+\cA)=& \bigl((a_{3}-a_{12})^2+(a_{2}+a_{13})^2+(a_{1}-a_{23})^2\bigr)\\
   &\quad\times\bigl(a_{3}+a_{12})^2+(a_{2}-a_{13})^2+(a_{1}+a_{23})^2\bigr)
   = a_{+}^2  a_{-}^2,
\end{split}
\end{equation}
from which follows that special cases will arise when
$\Det(\ba+\cA)=0$. Since the formulas for $a_{+}$ and $a_{-}$ are
expressed through  square roots, it is interesting to find a MV to
which the square roots are  associated. In
references~\cite{AcusDargysPreprint2020,Dargys2019} an algorithm
to compute the square root of MV in 3D algebras is provided. It
seems reasonable to conjecture that the special cases in
exponential are related to  isolated square roots  of the center
$a_S+a_I I$ of the considered algebra, where the scalars $a_S$ and
$a_I$ are defined by
\begin{equation}\begin{split}\label{aSaI}
&a_S=-(\ba+\cA)\d(\ba+\cA)=a_{1}^2+a_{2}^2+a_{3}^2+a_{12}^2+a_{13}^2+a_{23}^2,\\
&a_I =-(\ba+\cA)\w(\ba+\cA) I= -2(a_{3} a_{12}- a_{2} a_{13}+
a_{1} a_{23}).
\end{split}\end{equation}
In $\cl{0}{3}$ algebra the explicit formula for the center is
$-(\ba+\cA)(\ba+\cA)=a_S+a_I I$. In particular, the square root of
the center can be written as
\begin{equation}
\label{sSsqrtExprGenCl03}
  \begin{split}
\sqrt{a_S+a_I I}=&a_R+a_P I,\quad\textrm{where}\\
a_R+a_P I=&\begin{cases}
  \pm\frac{a_S+\sqrt{a_S^2-a_I^2}+a_I I}{\sqrt{2} \sqrt{a_S+\sqrt{a_S^2-a_I^2}}},\\[13pt]
  \pm\frac{a_S-\sqrt{a_S^2-a_I^2}+a_I I}{\sqrt{2} \sqrt{a_S-\sqrt{a_S^2-a_I^2}}},
\end{cases}\quad\text{if }  a_S^2>a_I^2\,.
  \end{split}
\end{equation}
From this follows that  $a_{+}$ and $a_{-}$ in Eq.~\eqref{exp03F}
can be expressed as $a_{+}^2=a_S+a_I= (a_R+a_P)^2$ and
$a_{-}^2=a_S-a_I=(a_R-a_P)^2$. Note that in
\eqref{sSsqrtExprGenCl03} the both required conditions $a_S>0$ and
$a_S^2>a_I^2$ are satisfied for all values of MV coefficients,
except when the vector and bivector parts of MV are absent. From
this we conclude that the condition $a_{+}=0$, or $a_{-}=0$, is
 equivalent to the determinant being zero, $\Det (a_S + a_I I)=(a_S + a_I)^2 (a_S- a_I)^2 = a_{+}^4  a_{-}^4=0$.

 The special cases in Theorem~\ref{exp03F}
occur when whichever of denominators, $a_{+}$ or $a_{- }$,
in the coefficients  turns to zero. Though at first glance we
could compute corresponding limits, for example, $\lim_{a_{+}\to
0}\frac{\sin a_{+}}{a_{+}}=1$,  in fact, the formula in this case
becomes  simpler because the condition $a_{+}=0$ implies that
$a_{3}=a_{12}, a_{2}=-a_{13}$ and $a_{1}=a_{23}$. Therefore, the
terms in vector and bivector components that include corresponding
differences vanish altogether.  Similarly, the case $a_{-}=0$
implies three conditions $a_{3}=-a_{12}, a_{2}=a_{13}$ and
$a_{1}=-a_{23}$ that nullify the corresponding terms in vector and
bivector components too. On the other hand, in scalar and
pseudoscalar components we can simply replace corresponding $\cos
a_{+}$ and $\cos a_{-}$ by $1$. Thus, the listed special cases
actually  represent the special cases already  found in the
analysis of algorithm of MV square root
in~\cite{AcusDargysPreprint2020}. After identification of
$a_0$ and $a_{123}$ with coefficients
in~\cite{AcusDargysPreprint2020}, $a_0\equiv s$ and $a_{123}\equiv
S$,  we find the following equivalence relations
$a_{+}^2=a_S+a_I=0\Leftrightarrow s=-S$,
$a_{-}=a_S-a_I=0\Leftrightarrow s=S$ and
$a_{-}=a_{+}=0\Leftrightarrow s=S=0$, respectively.

\section{MV exponentials in \cl{3}{0} and \cl{1}{2} algebras}
\label{sec:Cl30}
\begin{thm}[Exponential function in \cl{3}{0} (upper) and \cl{1}{2} (lower signs)]
\label{exp30thm} The exponential of MV
$\m{A}=a_0+a_1\e{1}+a_2\e{2}+a_3\e{3}+a_{12}\e{12}+a_{13}\e{13}+a_{23}\e{23}+a_{123}
I$ is another MV
$\exp(\m{A})=\ee^{a_{0}}\Bigl(b_0+\frac{b_1}{|c|}\e{1}+\frac{b_2}{|c|}\e{2}+\frac{b_3}{|c|}\e{3}+\frac{b_{12}}{|c|}\e{12}+\frac{b_{13}}{|c|}\e{13}+\frac{b_{23}}{|c|}\e{23}+b_{123}I\Bigr)$,
  where real coefficients $b_{i\ldots j}$ are
\begin{align}\label{exp30F}
&\kern-1em\begin{aligned}
    \phantom{{}_{99}} b_{0}=&\cos a_{123} \cos a_{-} \cosh a_{+}-\sin a_{123} \sin a_{-} \sinh a_{+},&&\notag\\
    b_{123}=&\sin a_{123} \cos a_{-} \cosh a_{+}+\cos a_{123} \sin a_{-} \sinh a_{+},&&
\end{aligned}&\notag\allowdisplaybreaks\\[3pt]
  &\begin{aligned}
 b_{1}=&\cosh a_{+} \sin a_{-} \bigl((a_{-} a_{1}-a_{+} a_{23}) \cos a_{123}-(a_{+} a_{1}+a_{-} a_{23}) \sin a_{123}\bigr)
  \\ &+\sinh a_{+}\cos a_{-} \bigl((a_{+} a_{1}+a_{-} a_{23}) \cos a_{123}+(a_{-} a_{1}-a_{+} a_{23}) \sin a_{123}\bigr),&&\\
  b_{2}=&\pm\cosh a_{+} \sin a_{-} \bigl((\pm a_{-} a_{2}+a_{+} a_{13}) \cos a_{123}+(\mp a_{+} a_{2}+a_{-} a_{13}) \sin a_{123}\bigr)\\
  &+\sinh a_{+}\cos a_{-} \bigl((a_{+} a_{2}\mp a_{-} a_{13}) \cos a_{123}+(a_{-} a_{2}\pm a_{+} a_{13}) \sin a_{123} \bigr),&&\\
  b_{3}=&\cosh a_{+} \sin a_{-} \bigl((a_{-} a_{3}\mp a_{+} a_{12}) \cos a_{123}\mp(\pm a_{+} a_{3}+a_{-} a_{12}) \sin a_{123}\bigr)\\
  &+\sinh a_{+}\cos a_{-} \bigl((a_{+} a_{3}\pm a_{-} a_{12}) \cos a_{123}+(a_{-} a_{3}\mp a_{+} a_{12}) \sin a_{123}\bigr),&&\\
\end{aligned}&\notag\allowdisplaybreaks\\[3pt]
  &\begin{aligned}
b_{12}=&\cosh a_{+} \sin a_{-} \bigl((\pm a_{+} a_{3}+a_{-} a_{12}) \cos a_{123}\pm (a_{-} a_{3}\mp a_{+} a_{12}) \sin a_{123}\bigr)\\
  & +\sinh a_{+}\cos a_{-} \bigl((\mp a_{-} a_{3}+a_{+} a_{12}) \cos a_{123}+(\pm a_{+} a_{3}+a_{-} a_{12}) \sin a_{123}\bigr),&&\\
b_{13}=&\mp\cosh a_{+} \sin a_{-} \bigl((a_{+} a_{2}\mp a_{-} a_{13}) \cos a_{123}+(a_{-} a_{2}\pm a_{+} a_{13}) \sin a_{123}\bigr)\\
  & + \sinh a_{+}\cos a_{-} \bigl((\pm a_{-} a_{2}+a_{+} a_{13}) \cos a_{123}+(\mp a_{+} a_{2}+a_{-} a_{13}) \sin a_{123}\bigr),&&\\
  b_{23}=&\cosh a_{+} \sin a_{-} \bigl((a_{+} a_{1}+a_{-} a_{23}) \cos a_{123}+(a_{-} a_{1}-a_{+} a_{23}) \sin a_{123}\bigr)\\
  &+\sinh a_{+}\cos a_{-} \bigl((-a_{-} a_{1}+a_{+} a_{23}) \cos a_{123}+(a_{+} a_{1}+a_{-} a_{23}) \sin a_{123}\bigr),&&\\
\end{aligned}&\allowdisplaybreaks\\[3pt]
  \quad &\textrm{with }\notag\\
  &\begin{aligned}
    |c|=&\sqrt{a_S^2+a_I^2}=a_{+}^2+a_{-}^2;\quad\textrm{where}\quad  a_S=a_{1}^2\pm a_{2}^2\pm a_{3}^2\mp a_{12}^2\mp a_{13}^2-a_{23}^2;&&\notag\\
    &\hphantom{\sqrt{a_S^2+a_I^2}=a_{+}^2+a_{-}^2;\quad\textrm{where}}\quad a_I=2(a_{3} a_{12}- a_{2} a_{13}+ a_{1} a_{23}),\quad \textrm{and}&&\notag\\
\end{aligned}&\\
  &\begin{aligned}
  a_{+}=&\begin{cases}\tfrac{1}{\sqrt{2}}\sqrt{a_S+|c|}& a_I\neq 0\\
  \sqrt{a_S}& a_I= 0\textrm{ and } a_S>0\\
0& a_I= 0\textrm{ and } a_S<0
\end{cases}
\end{aligned}&\notag\\
  &\begin{aligned}
    a_{-}=&\begin{cases}\tfrac{1}{\sqrt{2}}\frac{a_I}{\sqrt{a_S+|c|}}\hphantom{|c|}& a_I\neq 0\\
 0& a_I= 0\textrm{ and } a_S>0\\
\sqrt{-a_S}& a_I= 0\textrm{ and } a_S<0
\end{cases}
\end{aligned}&
\label{aplisaminus}
 \end{align}

  When both $a_{+}=0$ and $a_{-}=0$, or alternatively both $a_{S}=0$ and
$a_{I}=0$, the formulas  are associated with special cases
considered below in the Subsec.~\ref{specialCL30CL13}.
\end{thm}
%
 \textit{Proof}\\
It is enough to check the defining property
Eq.~\ref{definingProperty}. The validity was also checked by
expanding in Taylor series up to order 6 with symbolic
coefficients and up to order 20 using random integers.

Since both \cl{1}{2} and \cl{3}{0} algebras are represented by
\bbC(2) matrices they are mutually isomorphic. Therefore, the same
formula may be used for $\cl{3}{0}$ and $\cl{1}{2}$ algebras
without modification if one takes into account one-to-one
equivalence. For example, either
\[\e{2}\leftrightarrow\e{13},\qquad\e{3}\leftrightarrow\e{12},\] or,
alternatively, \[\e{1}\leftrightarrow\e{12},\quad
\e{2}\leftrightarrow\e{13},\quad
\e{3}\leftrightarrow\e{1},\quad\e{12}\leftrightarrow\e{23},\quad\e{13}\leftrightarrow\e{2},
\quad\e{23}\leftrightarrow\e{3}.\] Those not explicitly listed
being the same.

\subsection{Special cases of Theorem~\ref{exp30thm}}\label{specialCL30CL13}
The determinant of sum of vector and bivector parts of
MV $\m{A}$ in this case is
\begin{equation}
 \begin{split}\label{aPaM30}
   \Det(\ba+\cA)=& \Bigl(4\bigl(a_{3} a_{12}- a_{2} a_{13}+
a_{1} a_{23}\bigr)^2\\
   &\quad +\bigl(a_{1}^2\pm a_{2}^2\pm a_{3}^2\mp a_{12}^2\mp a_{13}^2-a_{23}^2\bigr)^2\Bigr)^2\\
   =& \bigl(a_{S}^2 + a_{I}^2\bigr)^2= a_{+}^2 + a_{-}^2,
\end{split}
\end{equation}
where upper signs is for \cl{3}{0} and lower for \cl{1}{2}
algebra. Equation~\eqref{aplisaminus}  shows that special cases
occur again when $\Det(\ba+\cA)=0$. The isolated square roots of
$c=a_S+a_I I$ of \cl{3}{0} algebra are given by (both signs for
both algebras)
\begin{equation}
\label{sSsqrtExprGenCl30}
  \sqrt{c} =\sqrt{a_S+a_I I}=\pm\frac{a_S+\sqrt{a_S^2+a_I^2}+a_I I}{\sqrt{2} \sqrt{a_S+\sqrt{a_S^2+a_I^2}}} =\pm (a_{+}+a_{-}I),
\end{equation}
where the root  $\sqrt{a_S^2+a_I^2}$ is a norm:
$|c|=\sqrt{c\tilde{c}}=\sqrt{a_S^2+a_I^2}=a_{+}^2+a_{-}^2$. The
coefficients $a_S$ and $a_I$ represent coefficients at scalar and
pseudoscalar of geometric product $\ba+\cA$ by itself. In
particular, for $\cl{3}{0}$ algebra the explicit form is
$(\ba+\cA)(\ba+\cA)=a_S+a_I I$, where $a_S$ and $a_I$ are
expressed through  inner and outer products,
$a_S=(\ba+\cA)\d(\ba+\cA)=a_{1}^2\pm a_{2}^2\pm a_{3}^2\mp
a_{12}^2\mp a_{13}^2-a_{23}^2$ (upper signs for \cl{3}{0} and
lower for \cl{1}{2} algebra) and  $a_I =-(\ba+\cA)\w(\ba+\cA) I=
2(a_{3} a_{12}- a_{2} a_{13}+ a_{1} a_{23}) $.

The denominator in \eqref{exp30F} vanishes when
$|c|=\sqrt{a_S^2+a_I^2}=a_{+}^2+a_{-}^2=0$. It is easy to see that
in this case all vector and bivector coefficients become zero
$b_{1}=b_{2}=b_{3}=b_{12}=b_{13}=b_{23}=0$. Then, in the
expressions for $b_{0}$ and $b_{123}$, we have to take $\cosh
a_{+} =\cos a_{-}=1$ and $\sinh a_{+} =\sin a_{-}=0$. After
identification with  coefficients
of~\cite{AcusDargysPreprint2020}, $a_0\equiv s$ and $a_{123}\equiv
S$, these conditions again are analogues of the only possible
special case when  $s=S=0$ in the square root of MV for
\cl{3}{0}~\cite{AcusDargysPreprint2020}.

\section{MV exponential in \cl{2}{1} algebra}
\label{sec:Cl21}

\begin{thm}[Exponential function in \cl{2}{1}]
\label{exp21thm} Exponential of  MV
$\m{A}=a_0+a_1\e{1}+a_2\e{2}+a_3\e{3}+a_{12}\e{12}+a_{13}\e{13}+a_{23}\e{23}+a_{123}
I$ is another MV $\exp(\m{A})= b_0+b_1\e{1}+b_2\e{2}+b_3\e{3}+b_{12}\e{12}+b_{13}\e{13}+b_{23}\e{23}+b_{123}I$, where

\begin{align}\label{exp21F}
 &\begin{aligned}
 \phantom{{}_{99}}b_{0}&=\frac{1}{2}\ee^{a_{0}}\bigl(\ee^{a_{123}}\mathop{\mathrm{co}}(a_{+}^2) + \ee^{-a_{123}} \mathop{\mathrm{co}}(a_{-}^2)\bigr),\\
b_{123}&=\frac{1}{2}\ee^{a_{0}}\bigl(\ee^{a_{123}}\mathop{\mathrm{co}}(a_{+}^2) -\ee^{-a_{123}}\mathop{\mathrm{co}}(a_{-}^2)\bigr),
\end{aligned}\notag\allowdisplaybreaks\\[3pt]
&\begin{aligned}
\phantom{{}_{99}}b_{1}&=\frac{1}{2}\ee^{a_{0}}\Bigl(\ee^{a_{123}}(a_{1}+a_{23}) \mathop{\mathrm{si}}(a_{+}^2) + \ee^{-a_{123}}(a_{1}-a_{23}) \mathop{\mathrm{si}}(a_{-}^2)\Bigr),\\
\phantom{{}_{99}}b_{2}&=\frac{1}{2}\ee^{a_{0}}\Bigl(\ee^{a_{123}}(a_{2}-a_{13})\mathop{\mathrm{si}}(a_{+}^2) +\ee^{-a_{123}}(a_{2}+a_{13})\mathop{\mathrm{si}}(a_{-}^2)\Bigr),\\
\phantom{{}_{99}}b_{3}&=\frac{1}{2}\ee^{a_{0}}\Bigl(\ee^{a_{123}}(a_{3}-a_{12})\mathop{\mathrm{si}}(a_{+}^2) +\ee^{-a_{123}}(a_{3}+a_{12})\mathop{\mathrm{si}}(a_{-}^2)\Bigr),\end{aligned}\notag\allowdisplaybreaks\\[3pt]
  &\begin{aligned}
    \phantom{{}_{9}}b_{12}&=\frac{1}{2}\ee^{a_{0}}\Bigl(-\ee^{a_{123}}(a_{3}-a_{12})\mathop{\mathrm{si}}(a_{+}^2) +\ee^{-a_{123}}(a_{3}+a_{12})\mathop{\mathrm{si}}(a_{-}^2)\Bigr),\\
    \phantom{{}_{99}} b_{13}&=\frac{1}{2}\ee^{a_{0}}\Bigl(-\ee^{a_{123}}(a_{2}-a_{13})\mathop{\mathrm{si}}(a_{+}^2) +\ee^{-a_{123}}(a_{2}+a_{13})\mathop{\mathrm{si}}(a_{-}^2)\Bigr),\\
 \phantom{{}_{99}}   b_{23}&=\frac{1}{2}\ee^{a_{0}}\Bigl(\ee^{a_{123}}(a_{1}+a_{23})\mathop{\mathrm{si}}(a_{+}^2) -\ee^{-a_{123}}(a_{1}-a_{23})\mathop{\mathrm{si}}(a_{-}^2)\Bigr),
\end{aligned}\notag \allowdisplaybreaks\\
 \quad \textrm{with}\\
  a_{+}^2&= -(a_{3}-a_{12})^2+(a_{2}-a_{13})^2+(a_{1}+a_{23})^2,\notag \\
  a_{-}^2&=-(a_{3}+a_{12})^2+(a_{2}+a_{13})^2+(a_{1}-a_{23})^2\,, \quad \textrm{and}\notag \\[4pt]
  \mathop{\mathrm{si}}(a_{\pm}^2)=&\begin{cases}
    \frac{\sinh \sqrt{a_{\pm}^2}}{\sqrt{a_{\pm}^2}},& a_{\pm}^2 > 0,\\
    \frac{\sin \sqrt{-a_{\pm}^2}}{\sqrt{-a_{\pm}^2}},& a_{\pm}^2 <0, \\
  \end{cases}\qquad \mathop{\mathrm{co}}(a_{\pm}^2)=\begin{cases}
   \cosh \sqrt{a_{\pm}^2},& a_{\pm}^2 > 0,\\
   \cos \sqrt{-a_{\pm}^2},& a_{\pm}^2 <0. \\
  \end{cases}
\end{align}
  When either $a_{+}^2=0$ or $a_{-}^2=0$ or both are zeroes,
  the formula yields special cases considered in Subsec~\ref{exp21thm}.
\end{thm}
\textit{Proof} The same as for $\cl{0}{3}$ and $\cl{3}{0}$
algebras, see Eq.~\eqref{definingProperty}.

\subsection{Special cases of Theorem~\ref{exp21thm}}

Determinant of the sum of vector and bivector of $\m{A}$ yields
\begin{equation}
 \begin{split}\label{aPaM21}
   \Det(\ba+\cA)=& \bigl(-(a_{3}-a_{12})^2+(a_{2}-a_{13})^2+(a_{1}+a_{23})^2\bigr)\\
   &\times\bigl(-(a_{3}+a_{12})^2+(a_{2}+a_{13})^2+(a_{1}-a_{23})^2\bigr)= a_{+}^2  a_{-}^2.
 \end{split}
\end{equation}
The special cases  occur when $\Det(\ba+\cA)=0$. As for previous
algebras they are related to the isolated roots of the element of
the center $a_S+a_I I$ of~\cl{2}{1}. In particular for the root of
the center  we find
\begin{equation}
\label{sSsqrtExprGenCl21}
  \begin{split}
\sqrt{a_S+a_I I}=&a_R+a_P I,\\
a_R+a_P I=&
    \begin{cases}
  \pm\frac{a_S+\sqrt{a_S^2-a_I^2}+a_I I}{\sqrt{2} \sqrt{a_S+\sqrt{a_S^2-a_I^2}}},\quad\text{if } a_S+\sqrt{a_S^2-a_I^2}>0\text{ and }
  a_S^2>a_I^2,
  \\[13pt]
  \pm\frac{a_S-\sqrt{a_S^2-a_I^2}+a_I I}{\sqrt{2} \sqrt{a_S-\sqrt{a_S^2-a_I^2}}},\quad\text{if } a_S-\sqrt{a_S^2-a_I^2}>0\text{ and }
  a_S^2>a_I^2.
\end{cases}
  \end{split}
\end{equation}
So, in $\cl{2}{1}$ algebra we have up to four roots. The real
coefficients $a_S$ and $a_I$ are equal to  coefficients  (which
are elements of the algebra center) of geometric product $\ba+\cA$
by itself. In particular, for $\cl{2}{1}$ algebra the explicit
form is $(\ba+\cA)(\ba+\cA)=a_S+a_I I$, where
$a_S=(\ba+\cA)\d(\ba+\cA)=a_{1}^2+a_{2}^2-a_{3}^2-a_{12}^2+a_{13}^2+a_{23}^2$
and  $a_I =(\ba+\cA)\w(\ba+\cA) I= 2(a_{3} a_{12}- a_{2} a_{13}+
a_{1} a_{23}) $.

In~\eqref{exp21thm}, $a_{+}$ and $a_{-}$ then again can be
expressed as $a_{+}^2=a_S+a_I= (a_R+a_P)^2=
-(a_{3}-a_{12})^2+(a_{2}-a_{13})^2+(a_{1}+a_{23})^2$ and
$a_{-}^2=a_S-a_I=
(a_R-a_P)^2=-(a_{3}+a_{12})^2+(a_{2}+a_{13})^2+(a_{1}-a_{23})^2$.
After comparison with $\cl{0}{3}$ algebra case, we see that the
explicit expressions now have different signs and, in general, can
acquire positive and negative values. Since these expressions are
present inside the square root  of  exponential, we formally have
to introduce functions $\mathop{\mathrm{si}}(a_{\pm}^2)$ and
$\mathop{\mathrm{co}}(a_{\pm}^2)$ (see Eq.~\ref{exp21thm})
in order to ensure real arguments for both functions.

When denominator $a_{+}$ or $a_{- }$ in Eq.~\eqref{exp21thm}
acquires zero value we have a special case. This corresponds to
the condition $\Det (a_S + a_I I)=(a_S - a_I)^2 (a_S + a_I)^2=0$.
Therefore,  conditions $a_{+}^2=a_S+a_I=0$ and $a_{-}^2=a_S-a_I=0$
define special cases. This requires to modify some of the  terms
in Eq.~\eqref{exp21thm}, i.e., these terms  have to be replaced by
limits $\lim_{a_{\pm}\to 0} \mathop{\mathrm{si}}(a_{\pm}^2)=1$.
Note, that now the coefficients  in vector and bivector components
that include $a_{+}$ or $a_{-}$, in general, do not necessary
vanish, unless the both $a_{+}^2$ and $a_{-}^2$ are equal to zero
simultaneously. This is a different situation compared to
$\cl{0}{3}$ algebra, for which the corresponding terms in the
component expressions always vanish.

Once more we stress that after identification of the coefficients
with those in~\cite{AcusDargysPreprint2020}, $a_0\equiv s$
and $a_{123}\equiv S$, the mentioned special cases correspond to
special cases of square root of MV, when
$a_{+}^2=a_S+a_I=0\Leftrightarrow s=-S$,
$a_{-}=a_S-a_I=0\Leftrightarrow s=S$ and
$a_{-}=a_{+}=0\Leftrightarrow s=S=0$, respectively.

\section{Particular cases: Pure bivector, vector and\\  (pseudo)scalar}
\label{sec:particularCases} In multivector $\m{A}$, equating
appropriate scalar coefficients  to zero, from formulas
\eqref{exp03F}, \eqref{exp30thm} and \eqref{exp21thm} one can
derive the exponentials of blades and compare them with those in
the literature, mainly for \cl{3}{0} and \cl{0}{3} algebras. For
mixed signature algebras, to authors knowledge, the formulas are
absent.

\subsection{Exponential of bivector}\label{expbivector}
In this case the exponential of a pure bivector $\cA=a_{12}
\e{12}+a_{13} \e{13}+a_{23} \e{23}$ can be expressed in a
coordinate-free form. The general formulas \eqref{exp03F},
\eqref{exp30thm} and \eqref{exp21thm} then reduce to
\begin{equation}\label{pureBivector}
\begin{split}
\ee^{\cA}=&\begin{cases}
\cos|\cA|+\frac{\cA}{|\cA|}\sin|\cA|&\textrm{for\ }\cl{3}{0},\cl{0}{3},\quad\cA^2<0, \\
\cosh|\cA|+\frac{\cA}{|\cA|}\sinh|\cA|&\textrm{for\ }\cl{1}{2},\cl{2}{1},\quad\cA^2>0,\\
\cos|\cA|+\frac{\cA}{|\cA|}\sin|\cA|&\textrm{for\ }\cl{1}{2},\cl{2}{1},\quad\cA^2<0,\\
\end{cases}
\end{split}
\end{equation}
where\  $|\cA|=\begin{cases}
  \sqrt{\cA^2}& \textrm{if}\quad\cA^2>0,\\
   \sqrt{-\cA^2}&\textrm{if}\quad \cA^2<0,\\
\end{cases}$\quad\textrm{and}
\begin{equation*}\cA^2=\begin{cases}
&-a_{12}^2-a_{13}^2-a_{23}^2\quad \textrm{\ for\ }\cl{3}{0},\cl{0}{3},\\
&-a_{12}^2+a_{13}^2+a_{23}^2\quad \textrm{\ for\ }\cl{2}{1},\\
&+a_{12}^2+a_{13}^2-a_{23}^2\quad \textrm{\ for\ }\cl{1}{2}.\\
 \end{cases}
\end{equation*}

\subsection{Exponential of vector}

In the case of pure vector $\ba=a_1\e{1}+a_2\e{2}+a_3\e{3}$, its
magnitude is $|\ba|=\sqrt{\pm\ba^2}$ where the root must be a
positive real number. Then, the general formulas reduce~to
\begin{equation}
\label{pureVector}
\begin{split}
\ee^{\ba}=&\begin{cases}
\cosh|\ba|+\frac{\ba}{|\ba|}\sinh|\ba|,&\quad \textrm{for}\quad \cl{3}{0}, \quad \ba^2>0,\\
\cos |\ba|+\frac{\ba}{|\ba|}\sin|\ba|,&\quad \textrm{for}\quad \cl{0}{3}, \quad \ba^2<0,\\
\cosh|\ba|+\frac{\ba}{|\ba|}\sinh|\ba|,&\quad \textrm{for}\quad\cl{1}{2}, \cl{2}{1},\quad \ba^2>0,\\
\cos|\ba|+\frac{\ba}{|\ba|}\sin|\ba|,&\quad \textrm{for}\quad\cl{1}{2}, \cl{2}{1}, \quad\ba^2<0,\\
\end{cases}
 \end{split}
\end{equation}
where $|\ba|=\begin{cases}
  \sqrt{\ba^2}& \ba^2>0,\\
   \sqrt{-\ba^2}& \ba^2<0,\\
\end{cases}$\quad and
\begin{equation*}
\ba^2= \begin{cases}
\pm(a_{1}^2+a_{2}^2+a_{3}^2),&\textrm{for}\quad \cl{3}{0}\textrm{\ ($+$ sign)\ },\cl{0}{3}\textrm{\ ($-$ sign)\ },\\
a_{1}^2+a_{2}^2-a_{3}^2,&\textrm{for}\quad \cl{2}{1},\\
a_{1}^2-a_{2}^2-a_{3}^2,&\textrm{for}\quad \cl{1}{2}.\\
\end{cases}
\end{equation*}
Thus, $|\ba|=\sqrt{a_1^2+a_2^2+a_3^2}$ for both \cl{3}{0} and
\cl{0}{3}.

\subsection{Exponent of scalar + pseudoscalar}
When $\m{A}=a_0+a_{123}I$, the type of the function depends on
sign of $I^2$, minus for \cl{3}{0} and \cl{1}{2}, and plus sign
for \cl{0}{3} and \cl{2}{1},
\begin{equation}
\label{scalarPseudosacalar}
\begin{split}
\ee^{a_0+a_{123}I}= &\begin{cases}
   \ee^{a_0} \bigl(\cos a_{123}+I \sin a_{123} \bigr),\qquad \textrm{for}\quad \cl{3}{0},\cl{1}{2},\\
   \ee^{a_0} \bigl(\cosh a_{123}+I \sinh a_{123}\bigr),\quad \textrm{for}\quad \cl{0}{3},\cl{2}{1}.\\
\end{cases}
\end{split}
\end{equation}
All listed in this section formulas are
well-known~\cite{Lounesto97}, and they readily follow from general
formulas \eqref{exp03F}, \eqref{exp30thm} and \eqref{exp21thm}.
One also can check that the identity
$\ee^{\m{A}}\ee^{-\m{A}}=\ee^{-\m{A}}\ee^{\m{A}}=1$ holds, i.e.
the inverse of exponential can be  obtained by changing the sign
of the exponent.

\section{Relations of the exponential to GA trigonometric and hyperbolic functions}
\label{Relations}

The geometric product is non-commutative. However, any two GA
functions of the same argument, for example $f(\m{A})$ and
$g(\m{A})$, that  can be expanded in the Taylor series, commute:
$f(\m{A})g(\m{A})=g(\m{A})f(\m{A})$. Indeed, for any chosen finite
series expansion  we have a product of two polynomials of a single
variable~$\m{A}$. Since the MV always commutes with itself, it
follows, that a well behaved functions of the same MV argument
commute too.

As known, the elementary trigonometric and hyperbolic  functions
in GA are defined by exactly the same series expansions as their
commutative counterparts
\cite{Lounesto97,Chappell2015,Josipovic2019}. GA hyperbolic
functions can be defined for an arbitrary MV  similarly as
ordinary functions, however, the GA trigonometric functions, in
general, only exist for real GAs that are characterized by a
commutative pseudoscalar  and property
$I^2=-1$~\cite{Chappell2015}, i.e., only for Clifford algebras
\cl{3}{0} and \cl{1}{2}. In order to define them for the algebras
\cl{0}{3} and \cl{2}{1} we have to introduce imaginary unit, i.e.
in these algebras trigonometric functions exist only when they are
complexified.

As known, scalar trigonometric and hyperbolic functions are linked
up through  the  imaginary unit $\ii=\sqrt{-1}$, for example
$\cos(\ii x)=\cosh(x)$ for all $x\in\bbR$. For MV functions
similar relations also exist if apart from $\ii$ the pseudoscalar
$I$ is included:
\begin{equation}
\cosh(I\m{A}) =\cos(\ii I\m{A}),\quad\sinh(I\m{A})=-\ii\sin(\ii
I\m{A}).
\end{equation}
Also, trigonometric and hyperbolic functions of MV $\m{A}$ can be
expressed through the exponentials, if one remembers that $I^2=-1$
for \cl{3}{0} and \cl{1}{2}, and $I^2=+1$ for \cl{0}{3} and
\cl{2}{1},
\begin{equation}\begin{split}\label{sincos}
& \sin\m{A}=\begin{cases}
   &\frac{I}{2}(\ee^{-I\m{A}}-\ee^{I\m{A}})\quad\text{\ for}\quad\cl{3}{0},\cl{1}{2},\\
   &\frac{\ii}{2}(\ee^{-\ii I\m{A}}-\ee^{\ii I\m{A}})\quad\text{for}\quad\cl{0}{3},\cl{2}{1},\\
   \end{cases}\\
& \cos\m{A}=\begin{cases}
   &\frac12(\ee^{-I\m{A}}+\ee^{I\m{A}})\quad\text{\ for}\quad\cl{3}{0},\cl{1}{2},\\
   &\frac12(\ee^{-\ii I\m{A}}+\ee^{\ii I\m{A}})\quad\text{for}\quad\cl{0}{3},\cl{2}{1},\\
    \end{cases}\\
\end{split}
\end{equation}
where $I\m{A}$ is the dual to multivector $\m{A}$. As suggested at
the beginning of this section the hyperbolic GA functions do not
require imaginary unit, thus we have
 \begin{equation}\label{sinhCoshDef}
\sinh\m{A} =\frac{1}{2}(\ee^{\m{A}}-\ee^{-\m{A}}),\quad \cosh\m{A}
=\frac{1}{2}(\ee^{\m{A}}+\ee^{-\m{A}}).
\end{equation}
From the above formulas follows  various relations between GA
trigonometric and hyperbolic functions that are analogues of the
well-known scalar relations, for example, a few of them   are
given below:
\begin{equation}\begin{split}
&\cos^2\m{A}+\sin^2\m{A}=1,\quad
\cosh^2\m{A}-\sinh^2\m{A}=1,\\
&\sin(2\m{A})=2\sin\m{A}\cos\m{A}=2\cos\m{A}\sin\m{A},\\
&\cos(2\m{A})=\cos^2\m{A}-\sin^2\m{A}.
\end{split}\end{equation}
Also, it should be noted that GA sine and cosine functions as well
as hyperbolic GA sine and cosine functions commute:
$\sin\m{A}\cos\m{A}=\cos\m{A}\sin\m{A}$ and
$\sinh\m{A}\cosh\m{A}=\cosh\m{A}\sinh\m{A}$.

Apart from relations between the  exact hyperbolic sine-cosine
functions and the exponential given in Eq.~\eqref{sinhCoshDef} we
can write  an exact formula for hyperbolic tangent as well (see
the beginning in this Section),
\begin{equation}\label{tanhDef}
\tanh \m{A}=\sinh \m{A}\cosh^{-1}\m{A}=\cosh^{-1} \m{A} \sinh
\m{A}\, ,
\end{equation}
and likewise for $\coth \m{A}$ functions. After substitution of
exponential formulas~\eqref{sinhCoshDef} into the right hand side
of \eqref{tanhDef} we obtain general  $\tanh \m{A}$. However, as a
first step in deriving exact formula for $\tanh \m{A}$ at first
one must compute the exact inverse of hyperbolic cosine. How to
compute the inverse MV in case of general Clifford algebras is
described in~\cite{Acus2018,Shirokov2020a,Helmstetter2019}. For
this purpose the adjoint and determinant of MV may be needed,
\begin{equation}\label{invDef}
\m{A}^{-1}=\frac{\Adj(\m{A})}{\Det(\m{A})}, \qquad
\Adj(\m{A})\,\m{A}=\m{A}\,\Adj(\m{A})=\Det(\m{A}).
\end{equation}
Here Det is the  determinant of MV, which in 3D can be computed
with the help of involutions~\cite{Dadbeh2011,Acus2018}
\begin{equation}\label{det3}
  \Det(\m{A})= \m{A} \reverse{\m{A}} \gradeinverse{\m{A}} \gradeinverse{\reverse{\m{A}}},
\end{equation}
where $\reverse{\m{A}}$ denotes reverse MV, and
$\gradeinverse{\m{A}}$ is grade inverse of MV $\m{A}$. Although
the computation of inverse of general 3D MV is straightforward,
the resulting symbolic expression is too large to be presented
here. However,  in the Appendix we shall make  use of numerical
calculations by \textit{Mathematica} for this purpose.

\section{Applications}
\label{sec:applications}

\subsection{Time-dependent GA equation with a simple Hamiltonian}
The spinor evolution  under the action of magnetic field is
considered. The field (vector) is assumed to consist of two parts,
constant parallel to $\e{3}$ and rotating in  $\e{12}$ plane with
angular frequency $\omega$,
\begin{equation}\label{Magnfield}
\bB(t)=B_0\e{3}+B_1\big(\e{1}\cos(\omega t)+\sigma\e{2}\sin(\omega
t)\big)
\end{equation}
The sign number $\sigma$  determines the  rotation sense. When
$\sigma=-1$ the  field of amplitude $B_1$ is rotating clockwise
and when $\sigma=1$ anticlockwise.

The time-dependent  Pauli-Schr{\"o}dinger equation in the presence
of homogeneous $\bB(t)$ field  for a spinor $\psi$, which is the
MV of \cl{3}{0} algebra, is
\begin{equation}\label{SchrodEq}
\frac{\dd\psi}{\dd t}=\frac{1}{2}\gamma I\bB(t)\psi,
\end{equation}
where $\gamma$ is the gyromagnetic ratio. This GA equation can be
solved by rotating frame method (in physics it is called the
rotating wave approximation) if the following rotor
$S=\exp(-\sigma\e{12}\omega t/2)$ is applied to
Eq.~\eqref{SchrodEq}. Multiplying  from left by reverse of $S$ and
then differentiating  with respect to time, we find
\begin{equation}\label{Spsi}\begin{split}
\frac{\dd(\widetilde{S}\psi)}{\dd t}=\frac{\dd\widetilde{S}}{\dd
t}\psi+\widetilde{S}\frac{\dd\psi}{\dd t} & =\tfrac12\big(\sigma
\e{12}\omega \widetilde{S}\psi+\widetilde{S}\gamma
I\bB(t)\psi\big)\\
 &=\tfrac12\big(\sigma\e{12}\omega+\gamma
I\widetilde{S}\bB(t)S\big)(\widetilde{S}\psi).
\end{split}\end{equation}
 When $\sigma=\pm 1$,  the product
$\widetilde{S}\bB(t)S=B_0\e{3}+B_1\cos(\sigma\omega
t)(\e{1}\cos\omega t+\sigma\e{2}\sin\omega t)+B_1\sin(\sigma\omega
t)(\e{2}\cos\omega t+\sigma\e{1}\sin\omega t)$ reduces to
time-independent field $\widetilde{S}\bB(t)S=B_1\e{1}+B_0\e{3}$.
Therefore the GA differential equation becomes
\begin{equation}\label{SpsiA}
\frac{\dd(\widetilde{S}\psi)}{\dd
t}=\frac12\big(\sigma\e{12}\omega+\e{23}\omega_1+\e{12}\omega_0\big)(\widetilde{S}\psi),
\end{equation}
where $\omega_0=\gamma B_0$, $\omega_1=\gamma B_1$. Since
Eq.~\eqref{SpsiA} has a constant MV coefficient its solution is
the exponential function,
\begin{equation}\label{SpsiB}
(\widetilde{S}\psi)=\exp\big(\tfrac12(\sigma\e{12}\omega+\e{23}\omega_1+\e{12}\omega_0)\big)(\widetilde{S}\psi)_0.
\end{equation}
At $t=0$ the initial MV is $(\widetilde{S}\psi)_0=\psi(0)$.
Multiplying from left by $S$ and expanding the second exponential
according to Subsec.~\ref{expbivector}, finally, we have
\begin{equation}\label{psit}
\psi=\ee^{-\sigma\e{12}\omega t/2}\Big(\cos\frac{\alpha
t}{2}+\alpha^{-1}(\e{23}\omega_1+\e{12}\big(\omega_0+\sigma\omega)\big)\sin\frac{\alpha
t}{2}\Big)\psi(0),
\end{equation}
where
$\alpha=\tfrac12\big((\sigma\omega+\omega_0)^2+\omega_1^2\big)^{1/2}$.

\begin{figure}[t]
\centering
\includegraphics[width=4.4cm]{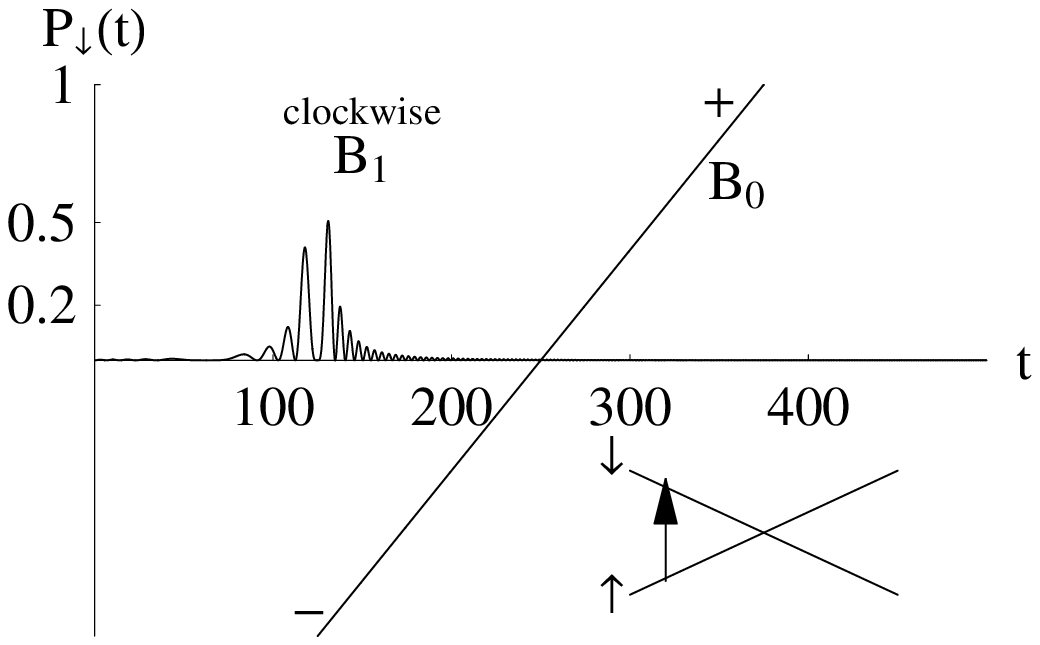}(a)\qquad
\includegraphics[height=2.8cm]{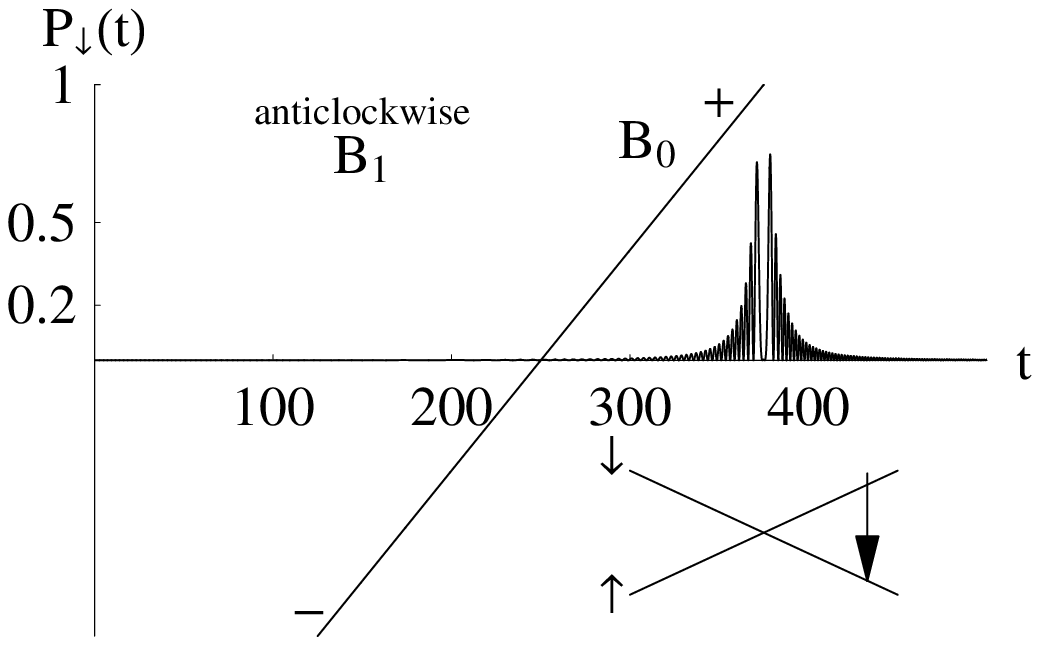}(b)
\caption{\footnotesize The probability $P_{\downarrow}(t)$ to find
the spin in the down direction when the magnetic field $\bB_0=B_0
\e{3}$ linearly increases from $B_0=-2$ to $B_0=2$ in the time
interval $T=(t_{fin}-t_{ini})=500$. The exciting field $\bB_1$
that flips the spin from $\uparrow$ to $\downarrow$ direction, as
shown in the insets  by long vertical arrows, is rotating in
$\e{12}$ plane, clockwise in~(a) and anticlockwise in (b). The
insets show the up and down spin eigenenergies  as a function of
magnetic field strength. In (a) the absorption and in (b) the
stimulated emission take place that are the main processes that
determine  performance of a laser which is the acronym of 'light
  absorption [and] stimulated emission [of] radiation'. Othor
parameters in the calculation: $\omega=1$, $\omega_1=0.05$.
\label{fig:spinTran}}
\end{figure}

Equation~\eqref{psit} describes the evolution of the total spinor
which is a mixture of up and down spinor  states
$\psi=\psi_{\uparrow}+\psi_{\downarrow}$ and normalized,
$\psi\widetilde{\psi}=1$. In GA the up and down spinor eigenstates
are, respectively, described by basis scalar $1$ and basis
bivector $\e{13}$ ~\cite{Doran03}. We shall assume that the spinor
initially is in the up eigenstate, $\psi(0)=\psi_{\uparrow}=1$.
Then,  the evolution of  the state $\psi_{\downarrow}$ is given by
projecting  $\psi$ onto the down eigenstate
$\e{13}$~\cite{Doran03}. The result is
\begin{equation}
\psi_{\downarrow}=-\langle\e{13}\psi\rangle+\langle\e{13}\psi\e{12}\rangle\e{12}=
\alpha^{-1}\Big(\sin\frac{\alpha t}{2}\big(\e{12}\cos\frac{\alpha
t}{2}-\sin\frac{\sigma\alpha t}{2}\big)\Big)
\end{equation}
The probability to detect the down spin at the moment $t$ then is
\begin{equation}
P_{\downarrow}(t)=\psi_{\downarrow}\widetilde{\psi}_{\downarrow}=
\left(\frac{\omega_1\sin\Big(\tfrac12
t\sqrt{(\sigma\omega+\omega_0)^2+\omega_1^2}\Big)}
{\sqrt{(\sigma\omega+\omega_0)^2+\omega_1^2}}\right)^2.
\end{equation}
At resonance, when $\sigma\omega+\omega_0=0$ (for clockwise
rotation $-\omega+\omega_0=0$, and for anticlockwise rotation
$\omega-\omega_0=0$) the probability oscillates,
$P_{\downarrow}(t)=\sin^2(\omega_1 t/2)=\sin^2(\gamma B_1 t/2)$,
with the frequency that depends on exciting field amplitude~$B_1$.
In quantum mechanics such rotating field induced oscillations
between up and down states  are called Rabi oscillations. If
magnetic field $\omega_0=\gamma B_0$ changes very slowly
(adiabatically) in the interval $T>>2\pi/\omega$, then in the
vicinity of resonance the probability peaks related with Rabi
oscillations will appear, Fig.~\ref{fig:spinTran}. The moment of
the appearance depends on the rotation sense via  sign
number~$\sigma$. The observed asymmetry between (a) and (b) panels
in  Fig.~\ref{fig:spinTran} is the manifestation of selection
rules for quantum transition under action  by rotating magnetic
field.

\subsection{Relations to even $\cl{1}{3}^{+}$ and $\cl{3}{1}^{+}$ geometric algebras}

If isomorphism rules between even subalgebra of $\cl{1}{3}$ and
full $\cl{3}{0}$ algebra are made of, namely,
\[\begin{array}{llll}
\cl{1}{3}^{+}\leftrightarrow\cl{3}{0},
&\e{23}\leftrightarrow\e{12},&\e{24}\leftrightarrow\e{13},
&\e{34}\leftrightarrow\e{23},\\ \e{1234}\leftrightarrow\e{123},
&\e{12}\leftrightarrow\e{1},&\e{13}\leftrightarrow\e{2},&\e{14}\leftrightarrow\e{3},.
\end{array}\]
or, alternatively,
\[\begin{array}{llll}
\cl{1}{3}^{+}\leftrightarrow\cl{3}{0}&\e{23}\leftrightarrow\e{23},&\e{24}\leftrightarrow\e{13},&
\e{34}\leftrightarrow\e{12},\\ \e{1234}\leftrightarrow\e{123},&
\e{12}\leftrightarrow\e{3},&
\e{13}\leftrightarrow\e{2},&\e{14}\leftrightarrow\e{1}.
\end{array}\] it is
easy to obtain explicit formulas for physically important cases of
exponentials of general even MVs that represents  spinors in
\cl{1}{3} algebra.

In case of \cl{3}{1} the following rules may be used for this
purpose,
\[\begin{array}{llll}
\cl{3}{1}^{+}\leftrightarrow\cl{3}{0},&\e{12}\leftrightarrow\e{12},&\e{13}\leftrightarrow\e{13},
&\e{23}\leftrightarrow\e{23},\\ \e{1234}\leftrightarrow\e{123},
&\e{14}\leftrightarrow\e{1},&\e{24}\leftrightarrow\e{2},
&\e{34}\leftrightarrow\e{3},\\
 \end{array}\]
or, alternatively,
\[\begin{array}{llll}
\cl{3}{1}^{+}\leftrightarrow\cl{3}{0},&\e{12}\leftrightarrow\e{23},&\e{13}\leftrightarrow\e{13},
&\e{23}\leftrightarrow\e{12},\\ \e{1234}\leftrightarrow\e{123},
&\e{14}\leftrightarrow\e{3},&
\e{24}\leftrightarrow\e{2},&\e{34}\leftrightarrow\e{1}.
  \end{array}\]

\section{Discussion and conclusions} \label{conclusion}

Since the obtained exponentials are expressed in coordinates the
final formulas appear rather complicated. In geometric Clifford
algebra the formulas in coordinate-free form may be desirable. The
main problem is with vectors and bivectors the components of
which, as seen from Eqs,~\eqref{exp03F}, \eqref{exp30thm} and
\eqref{exp21thm}, are entangled mutually. To avoid the
entanglement, a better strategy\footnote{It should be noted that
at present the existing symbolic packages can do calculations in a
concrete orthogonal frame (basis) rather than with simple blades
directly.} would be to avoid MV expansion in components  at all as
 done in~\cite{Chappell2015}.

Let's take \cl{3}{0} and introduce the following complex
quantity~\cite{Dargys}
\begin{equation}\begin{split}
&H=\Big((a_1-\ii a_{23})^2+(a_2-\ii a_{31})^2+(a_3-\ii a_{12})^2\Big)^{1/2}, \\
&H^*=\Big((a_1+\ii a_{23})^2+(a_2+\ii a_{31})^2+(a_3+\ii
a_{12})^2\Big)^{1/2},
\end{split}\end{equation}
so that $HH^*$ is a real number. Introduction of the imaginary
unit makes the formulas more compact and permits
trigonometric-hyperbolic expansion  of~$\ee^{\m{A}}$. In the
expanded form the vector and bivector coefficients in $HH^*$
represent the sum of $81$ terms that consist of various products
of $a_i$ and $a_{ij}$. However, the function $H$ can be written
very compactly if coordinate-free form is used~\cite{Dargys},
\begin{equation}
H=\ba^2+\cA^2+2\ii I\,\ba\w\cA.
\end{equation}

The same motive is seen in the coefficients $a_{+}$ and $a_{-}$
that appear in the theorem~\ref{exp03thm}. Furthermore the
coefficients  may be given a similar shape:
\begin{equation}
a_{+}=\sqrt{\ba^2-\cA^2+2I\ba\w\cA}\,,\qquad
 a_{-}=\sqrt{\ba^2-\cA^2-2I\ba\w\cA}\,.
\end{equation}
So, there appears a chance to construct a MV exponential functions
having a compact and coordinate-free forms which will be more
useful and efficient in various practical GA applications.

In conclusion, we have been able to expand  the GA  exponential
function of a general argument into MV in the coordinate form for
all four 3D Clifford geometric algebras. The expansion has been
applied to get exact expressions for trigonometric and hyperbolic
GA  functions and to investigate the convergence of respective
series. It was found that both trigonometric and hyperbolic GA
sine-cosine series convergence is satisfactory if GA series is
limited to more than 6 terms. However the convergence of tangent
series is slower, about 40 significant figures are needed to reach
similar precision. We think that such an expansion of the
exponential will be useful in solving GA differential
equations~\cite{Snygg12,Dargys17,Dargys09bCL}, in signal and image
processing, in automatic control and robotics~\cite{Lavor2018}.

\bibliographystyle{REPORT}

%
\bibliography{expdim3}

\appendix
\section{Numerical comparison between  exact formulas and their series expansion}

In this Appendix a comparison between  exact MV formulas obtained
in Sec.~{\ref{Relations}} and  finite series expansion is made.
The numerical form of MVs is used for this purpose. The knowledge
of exact formulas allows to investigate the rate of convergence of
finite GA trigonometric and hyperbolic series in \cl{3}{0}
algebra. The following MV
\begin{equation}\label{mvAp}
  \m{A}^\prime=\frac{1}{N}\bigr(4+\e{1}+3 \e{2}-5 \e{3}+10 \e{12}+9 \e{13}-9 \e{23}-4 I\bigl),
\quad I=\e{123}
\end{equation}
is used for this purpose where the integer numbers were generated
randomly. The  normalization factor $N$ helps to make
trigonometric series convergent. Up to 8~significant figures are
presented in numerical  evaluation of symbolic (exact) formulas
from Sec.~{\ref{Relations}}. Of course, obtained exact formulas
can be used to compute trigonometric functions of any MV, even if
respective Taylor series does not converge, for example, at
large coefficients and  $N=1$. Our primary intention here
is however to compare answers provided by exact formula and Taylor
series expansion.

\subsection{GA hyperbolic functions}

The trigonometric function series can be made to converge if in
\eqref{mvAp} we chose large enough $N$ but not too large.
We have found that the optimal factor must be larger than
the determinant norm of MV in Eq.~\ref{det3}. The norm is defined
as the determinant of~$\m{A}$ raised to fractional power $1/k$,
where $k=2^{\lceil n/2\rceil}$, i.e. $|\m{A}|=
(\Det(\m{A}))^{1/k}$. This  norm can be interpreted as a number of
MVs $\m{A}$ in a MV product needed to define $\Det(\m{A})$. In our
case, $\Det(\m{A})$ in Eq.~\eqref{det3} consists of geometric
product of four MVs, therefore, for 3D algebras ($n=3$) we have
$k=2^{\lceil 3/2\rceil}=2^{2}=4$ and the determinant norm is
$|\m{A}|=\sqrt[4]{\Det(\m{A})}$. For the chosen MV $\m{A}^\prime$
we find $\Det (\m{A}^\prime)=71129$ and
$|\m{A}^\prime|=\sqrt[4]{71129}\approx 16.33$. Since the strict
analysis of convergence\footnote{If, instead, for example, we
divide the MV by the largest coefficient in the considered MV,
then we immediately would find that $\tanh \m{A}$ series fails to
converge.} of multivector series is outside the scope of this
article, we will divide the chosen MV by the nearest larger
integer $17>16.33$. Due to multiplicative property of the
determinant $\Det(\m{A}\m{A})=\Det(\m{A})\Det(\m{A})$, division by
any scalar that is larger than the determinant norm factor
$1/|\m{A}|$ ensures, that determinants of series terms make a
decreasing sequence, i.e. $|\Det(A)|>|\Det(\m{A}\m{A})|>\cdots
>|\Det(\m{A}\m{A}\cdots\m{A})|$, and, therefore, we may anticipate that
MV series will tend to converge, or at least will yield meaningful
answer. For GA series we have profited by  standard
exponential, trigonometric and hyperbolic
series~\cite{Abramowitz}. In particular, for $\tanh \m{A}$ we have
used $\tanh
\m{A}=\m{A}-\tfrac{1}{3}\m{A}^3+\tfrac{2}{15}\m{A}^5-\tfrac{17}{315}\m{A}^7+\tfrac{62}{2835}\m{A}^9+\cdots$.

To illustrate, let's compute hyperbolic functions $\sinh \m{A}$,
$\cosh \m{A}$,  $\cosh^{-1}\m{A}$ and\footnote{The function
$\cosh^{-1}(\m{A})$ can be calculated from series
$\cosh^{-1}(\m{A})=\sum_{n=0}^{\infty}\textrm{E}_n\frac{\m{A}^n}{n!}$,
where $\textrm{E}_n$ are the Euler coefficients, and the condition
$\cosh^{-1}(\m{A})\cosh(\m{A})=1$. In fact, the latter condition
gives the Euler numbers. In case of inverse trigonometric function
we have
$\cos^{-1}(\m{A})=\sum_{n=0}^{\infty}\textrm{E}_n\frac{\m{A}^{2n}}{2n!}$
and $\cos^{-1}(\m{A})\cos(\m{A})=1$. Similar relations exist for
hyperbolic and trigonometric tangent functions but now instead of
Euler numbers there appear Bernoulli numbers
$\textrm{B}_n$~\cite{Abramowitz}.} $\tanh \m{A}$ of normalized MV
argument $\m{A}^{\prime\prime}$,
\begin{equation}\label{Aprimeprime}
\m{A}^{\prime\prime}=\frac{1}{17}(4+\e{1}+3
\e{2}-5 \e{3}+10 \e{12}+9 \e{13}-9 \e{23}-4 I).
\end{equation}
Substituting  $\m{A}^{\prime\prime}$ into exact symbolic formulas
Eqs.~\eqref{sinhCoshDef} and \eqref{tanhDef} (where inverse MV is
computed using \eqref{invDef} and \eqref{det3}) and then
evaluating exact expressions numerically up to 8~significant
figures (the  last digit is exact) we obtain
\begin{align*}
   \sinh\m{A}^{\prime\prime}& =
 \begin{aligned}[t]
  &0.0806082\phantom{\e{13}} &\kern-1em -&0.0230640 \e{1}&\kern-1em +&0.0787983 \e{2}&\kern-1em -&0.1724390 \e{3}\\
     \kern-1em +&0.5504206 \e{12}&\kern-1em +&0.4830460 \e{13}&\kern-1em -&0.4666026 \e{23}&\kern-1em -&0.2082492 I\,,\\
\end{aligned}\notag\\
\cosh\m{A}^{\prime\prime} &=
 \begin{aligned}[t]
  &0.6039792 &\kern-1em -&0.1111834 \e{1}&\kern-1em -&0.0900922 \e{2}&\kern-1em +&0.0825265 \e{3}\\
  \kern-1em +&0.1730832 \e{12}&\kern-1em +&0.1354867 \e{13}&\kern-1em -&0.1084358 \e{23}&\kern-1em -&0.2939648 I\,,\\
 \end{aligned}\\
\tanh\m{A}^{\prime\prime} &=
  \begin{aligned}[t]
    & 0.6231177 &\kern-1em +& 0.3099294 \e{1}&\kern-1em +&0.4271905 \e{2}&\kern-1em -&0.5723737 \e{3}\\
  \kern-1em +&0.4466951 \e{12}&\kern-1em +&0.4439088 \e{13}&\kern-1em -&0.4997530 \e{23}&\kern-1em +&0.0547345 I\,.
  \end{aligned}\notag
\end{align*}
For comparison we provide answers obtained by finite  series
expansions,
\begin{align*}
\sinh_6\m{A}^{\prime\prime} =&\begin{aligned}[t]
  &0.0806569&\kern-1em -&0.0229633 \e{1}&\kern-1em +&0.0788338 \e{2}&\kern-1em -&0.1724240 \e{3}\\
  \kern-1em +&0.5500202 \e{12}&\kern-1em +& 0.4827078 \e{13}&\kern-1em -& 0.4662941 \e{23}&\kern-1em -&0.2076350 I\,,
\end{aligned}\notag\\
\cosh_6\m{A}^{\prime\prime} =& \begin{aligned}[t]
  &0.6040721&\kern-1em -&0.1111303 \e{1}&\kern-1em -&0.0900394 \e{2}&\kern-1em +&0.0824681 \e{3}\\
  \kern-1em+&0.1730517 \e{12}&\kern-1em +&0.1354672 \e{13}&\kern-1em -&0.1084281 \e{23}&\kern-1em -&0.2939312 I\,,
\end{aligned}\notag\\
  \tanh_6\m{A}^{\prime\prime} =& \begin{aligned}[t]
    &0.7629316 &\kern-1em +&0.3616722 \e{1}&\kern-1em +&0.5029447 \e{2}&\kern-1em -&0.6765545 \e{3}\\
    \kern-1em +&0.5446755 \e{12}&\kern-1em +&0.5387139 \e{13}&\kern-1em -&0.6033886 \e{23}&\kern-1em -&0.1176009 I\,,
\end{aligned}\\
\tanh_{40}\m{A}^{\prime\prime} =& \begin{aligned}[t]
  &0.6231595 &\kern-1em +&0.3099902 \e{1}&\kern-1em +&0.4272145 \e{2}&\kern-1em -&0.5723697 \e{3}\\
  \kern-1em+&0.4464672 \e{12}&\kern-1em +&0.4437168 \e{13}&\kern-1em -&0.4995786 \e{23}&\kern-1em +&0.0550762
  I\,.
\end{aligned}\notag
\end{align*}
The subscripts at hyperbolic functions indicate the number of
terms that has been included in the summation of finite series  to
get the result. It can be  seen that $\tanh\m{A}$ converges much
slower than $\cosh\m{A}$ and $\sinh\m{A}$. The latters  are
directly related to exponential. For $\tanh\m{A}$  we have had to
include 50 terms to get six exact figures. If instead
in~\eqref{mvAp} we would take different factor
$N=\sqrt[4]{71129}$, and then try to compute $\tanh\m{A}$ by
by  standard (textbook) series expansion, then we would
immediately find that the series fails to converge, whereas exact
formula that follows from exponential  yields meaningful answer.
One can also easily check that all computed MV functions commute
pairwise up to the computed precision.

\subsection{GA trigonometric  functions}

In this Appendix we restrict ourselves to \cl{3}{0} algebra for
which $I^2=-1$. The  exact formulas  in the exponential form  for
$\sin \m{A}$ and $\cos \m{A}$ in Eqs.~\ref{sincos} have been used.
The numerical MV given by Eq.~\ref{Aprimeprime} was inserted to
find the following exact GA functions presented  below  with 8
significant figures,
\begin{align*}
   \sin\m{A}^{\prime\prime}& =
 \begin{aligned}[t]
  &0.4142215\phantom{\e{13}} &\kern-1em +&0.1561775 \e{1}&\kern-1em +&0.2887099 \e{2}&\kern-1em -&0.4312306 \e{3}\\
     \kern-1em +&0.6127064 \e{12}&\kern-1em +&0.5664210 \e{13}&\kern-1em -&0.5864014 \e{23}&\kern-1em -&0.2430952 I\,,\\
\end{aligned}\notag\allowdisplaybreaks\\
\cos\m{A}^{\prime\prime} &=
 \begin{aligned}[t]
  &1.3837580 &\kern-1em +&0.1075001  \e{1}&\kern-1em +&0.0726490 \e{2}&\kern-1em -&0.0516785 \e{3}\\
  \kern-1em -&0.2436586 \e{12}&\kern-1em -&0.1984718 \e{13}&\kern-1em +&0.1707105 \e{23}&\kern-1em +&0.4152926 I\,,\\
 \end{aligned}\allowdisplaybreaks\\
\tan\m{A}^{\prime\prime} &=
  \begin{aligned}[t]
    & 0.0520468 &\kern-1em -& 0.0321336 \e{1}&\kern-1em +&0.0568865 \e{2}&\kern-1em -&0.1373908 \e{3}\\
  \kern-1em +&0.4876809 \e{12}&\kern-1em +&0.4261388 \e{13}&\kern-1em -&0.4091069 \e{23}&\kern-1em -&0.1473168 I\,.
  \end{aligned}\notag
\end{align*}
On the other hand,  using series expansion of $\sin\m{A}$ ,
$\cos\m{A}$ and $\tan\m{A}$, we find,
\begin{align*}
   \sin_6\m{A}^{\prime\prime}& =
 \begin{aligned}[t]
  &0.4141938\phantom{\e{13}} &\kern-1em +&0.1560854 \e{1}&\kern-1em +&0.2886852 \e{2}&\kern-1em -&0.4312611 \e{3}\\
     \kern-1em +&0.6131181 \e{12}&\kern-1em +&0.5667705 \e{13}&\kern-1em -&0.5867229 \e{23}&\kern-1em -&0.2437297 I\,,\\
\end{aligned}\notag\allowdisplaybreaks\\
\cos_6\m{A}^{\prime\prime} &=
 \begin{aligned}[t]
  &1.3838520 &\kern-1em +&0.1075543  \e{1}&\kern-1em +&0.0727025 \e{2}&\kern-1em -&0.0517373 \e{3}\\
  \kern-1em -&0.2436926 \e{12}&\kern-1em -&0.1984933 \e{13}&\kern-1em +&0.1707199 \e{23}&\kern-1em +&0.4153298 I\,,\\
 \end{aligned}\allowdisplaybreaks\\
\tan_6\m{A}^{\prime\prime} &=
  \begin{aligned}[t]
    & 0.0958579 &\kern-1em +& 0.0035747 \e{1}&\kern-1em +&0.0832419 \e{2}&\kern-1em -&0.1588803 \e{3}\\
  \kern-1em +&0.4184797 \e{12}&\kern-1em +&0.3705885 \e{13}&\kern-1em -&0.3625310 \e{23}&\kern-1em -&0.0454115 I\,.
  \end{aligned}\notag\allowdisplaybreaks\\
  \tan_{40}\m{A}^{\prime\prime} &=
  \begin{aligned}[t]
    & 0.0522097 &\kern-1em -& 0.0320415 \e{1}&\kern-1em +&0.0569781 \e{2}&\kern-1em -&0.1374922 \e{3}\\
  \kern-1em +&0.4876273 \e{12}&\kern-1em +&0.4261060 \e{13}&\kern-1em -&0.4090946 \e{23}&\kern-1em -&0.1472580 I\,.
  \end{aligned}\notag
\end{align*}
The subscripts at trigonometric functions show the number of terms
that has been used in series expansion to get the result.

\end{document}